\documentclass[letterpaper, 11pt]{article}

\usepackage[margin=1in]{geometry}

\usepackage{url}
\usepackage{amssymb, amsmath, mathtools, amsfonts, amsthm, bm, dsfont}
\usepackage{float, xspace, graphicx, enumitem}
\usepackage{placeins,wrapfig}
\usepackage[normalem]{ulem}
\usepackage{booktabs, multirow, array}
\usepackage{algorithm, algpseudocode}
\usepackage{physics}

\usepackage{tikz}
\definecolor{darkgreen}{RGB}{0,140,0}
\usetikzlibrary{arrows.meta, matrix, arrows, calc, shapes, backgrounds}
\usetikzlibrary{shapes.callouts, decorations.text, decorations.markings, shapes.misc}
\usepackage{pgfplots}
\pgfplotsset{compat=1.17}

\usepackage[
    bookmarks,
    bookmarksdepth=subsubsection,
    colorlinks=true,
    allcolors=black,
    citecolor=red,
    linkcolor=blue,
    plainpages=false
]{hyperref}
\usepackage[capitalize]{cleveref}

\theoremstyle{plain}
\newtheorem{theorem}{Theorem}
\newtheorem{lemma}{Lemma}
\newtheorem{proposition}{Proposition}
\newtheorem{corollary}{Corollary}

\theoremstyle{definition}

\newtheorem{remark}{Remark}

\newif\ifhideproofs
\ifhideproofs
    \usepackage{environ}
    \NewEnviron{hide}{}

\fi

\newcommand{\Expect}{\mathbb{E}}

\newcommand{\Var}{\mathrm{Var}}
\newcommand{\diff}{\mathrm{d}}
\renewcommand{\norm}[1]{\left\lVert #1 \right\rVert}

\renewcommand{\bar}{\overline}

\newcommand{\xbf}{\mathbf{x}}
\newcommand{\ybf}{\mathbf{y}}

\title{Service location flexibility in dynamic vehicle routing}   
\date{}
\author{%
    Gérard P. Cachon, Bora Cetin, and Sophie H.~Yu%
    \thanks{All authors are with the Operations, Information and Decisions Department,
    The Wharton School, University of Pennsylvania, Philadelphia, PA 19104, USA.
    \texttt{\detokenize{{cachon,bcet78,hysophie}@wharton.upenn.edu}}.}%
}

\hypersetup{pdfauthor={Gérard P. Cachon, Bora Cetin, and Sophie H. Yu}}
\newcommand{\E}[1]{\Expect\!\left[#1\right]}
\newcommand{\Pp}{\mathbb P}
\newcommand{\FR}{\mathbf{FR}}

\begin{document}
\maketitle
\begin{abstract}
Parcel lockers and delivery robots let a vehicle serve customers at nearby locations. We study this \emph{service location flexibility} in the dynamic traveling repairman problem, where a single vehicle may serve each demand anywhere within a radius $r$ of its location. Arrivals form a Poisson process, demand locations are independent and uniform in the unit square, and on-site service times are i.i.d. with uniformly bounded second moments, where the radius $r$ may vary across systems as the on-site load $\rho$ tends to one. In this heavy-traffic regime, the optimal mean system time is of order $(1-\rho)^{-1}(1+\min\{1/r,1/(1-\rho)\})$, compared with $(1-\rho)^{-2}$ without service location flexibility \cite{BvR1991, Pavone2011}. Thus, flexibility reduces the order of optimal mean system time when $r$ stays fixed or shrinks more slowly than $1-\rho$. Moreover, when $r$ is on the order of $1-\rho$ or larger, the vehicle can attain the optimal order by repeating a fixed covering route that passes within $r$ of every location. For example, the vehicle can follow a fixed grid route through parcel lockers, with customers in each square cell sharing the locker at its center.

\end{abstract}


\setcounter{tocdepth}{2}
\section{Introduction}

When two neighbors choose the same shop for parcel pickup, a driver can make both deliveries in one stop, though each parcel still takes time to handle. In a 2026 DHL survey, 29\% of respondents worldwide and 40\% of those in Europe said they usually received their online orders at a locker or parcel shop~\cite{dhlecommerce2026ooh}. We model this service location flexibility by allowing the driver to deliver each parcel anywhere within an agreed distance of the recipient's address.

For a given list of deliveries, a carrier can measure the benefit of this flexibility in miles saved \cite{Krushinsky2021}. When orders arrive over time, serving two nearby customers at a shared stop can shorten their wait and let the driver reach other waiting customers sooner. As the delivery system becomes more congested, how much can service location flexibility shorten the average time from order arrival to delivery?

We study this question in the dynamic traveling repairman model of Bertsimas and van Ryzin \cite{BvR1991}. Demands arrive at rate $\lambda$ according to a Poisson process, with independent uniform locations in the unit square. A single vehicle travels at unit speed and serves one demand at a time. The demands have independent, identically distributed on-site service times with mean $\bar s$. The service radius $r$ measures the extent of service location flexibility and is chosen before operation. Each demand must be served within distance $r$ of its address. The dispatcher chooses when to give the vehicle its next assignment and which waiting demands to include. We measure system time from a demand's arrival until the vehicle completes service at the chosen service point.

Consider two extremes as the on-site load $\rho=\lambda\bar s$ tends to one. If $r=0$, the vehicle must visit every address, and the optimal mean system time is of order $(1-\rho)^{-2}$~\cite{BvR1991,Pavone2011}. If $r\ge1/\sqrt{2}$, every customer can reach the center of the square. The vehicle can stay there and serve demands in arrival order, giving an M/G/1 queue with mean system time of order $(1-\rho)^{-1}$~\cite{BertsekasGallager1992}. As the load rises, more orders await delivery at the same time. The driver then has more opportunities to find customers who live within a short walk of the same stop. We therefore study how the reduction in mean system time depends on the service radius, allowing both constant radii and radii that shrink as the load rises.

In our paper, we characterize this transition for every service radius $r>0$. In heavy traffic, the optimal mean system time among all admissible stable policies has order
$\tfrac{1}{1-\rho}\bigl(1+\min\{\tfrac{1}{r},\tfrac{1}{1-\rho}\}\bigr)$.
We provide policies that attain this order and prove a matching lower bound for every admissible stable policy, including those that choose their routes dynamically based on the locations of waiting demands.
Compared to the $(1-\rho)^{-2}$ benchmark with no service location flexibility~\cite{BvR1991,Pavone2011}, the optimal mean system time is of smaller order with service location flexibility precisely when the service radius satisfies $r\gg1-\rho$.

In particular, we show that \emph{fixed-route policies} can attain the optimal order of mean system time when $r$ is of order $1-\rho$ or larger. Before demands arrive, the carrier chooses an \emph{$r$-covering route}, a continuous closed curve of finite length $L_r>0$ in the square that passes within distance $r$ of every address. On this route the carrier designates a depot and assigns each address a service point within distance $r$. Each cycle starts and ends at the depot. The vehicle traverses the full route and serves all demands waiting at departure. Demands arriving during the cycle are served on the next one, which begins immediately. For a fixed direction or an independent fair coin at each departure, we give the exact mean system time at every load $0<\rho<1$. If the vehicle chooses its direction by an independent fair coin toss at each departure, mean system time exceeds the non-spatial M/G/1 benchmark~\cite[Section~3.5]{BertsekasGallager1992} by exactly $L_r/(1-\rho)$. This excess is the mean duration of a cycle from the depot back to the depot, including travel and all on-site service. Always traveling in the same direction changes mean system time by at most $L_r/2$ relative to the fair-coin rule.

\begingroup
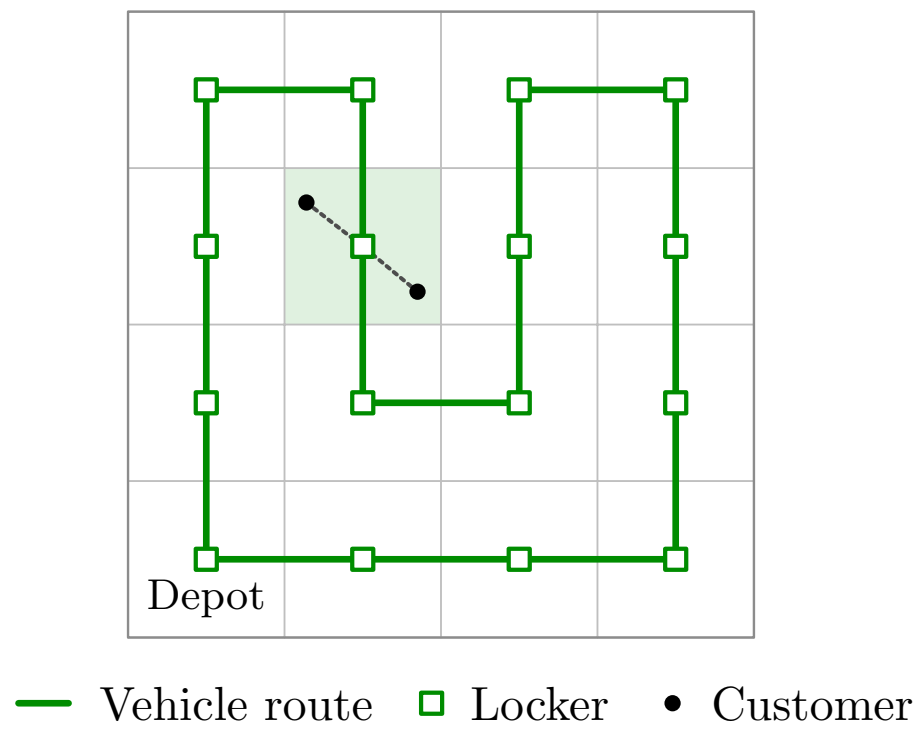
\begin{wrapfigure}{r}{0.45\textwidth}
\small
\centering
\begin{tikzpicture}[
    scale=0.9,transform shape,
    x=1.15cm,y=1.15cm,font=\small,line join=round,line cap=round,
    route/.style={draw=darkgreen,line width=1.25pt},
    locker/.style={rectangle,draw=darkgreen,fill=white,line width=0.8pt,
                   minimum size=4.4pt,inner sep=0pt},
    customer/.style={circle,fill=black,minimum size=3.4pt,inner sep=0pt}
]
    \fill[darkgreen!12] (1,2) rectangle (2,3);
    \draw[xstep=1.15cm,ystep=1.15cm,black!25,line width=0.35pt] (0,0) grid (4,4);
    \draw[black!45,line width=0.45pt] (0,0) rectangle (4,4);

    \draw[route] (0.5,0.5) -- (0.5,3.5) -- (1.5,3.5) -- (1.5,1.5)
        -- (2.5,1.5) -- (2.5,3.5) -- (3.5,3.5) -- (3.5,0.5) -- cycle;
    \coordinate (shared) at (1.5,2.5);
    \coordinate (first) at (1.14,2.78);
    \coordinate (second) at (1.85,2.21);
    \draw[black!70,densely dotted,line width=0.75pt]
        (first) -- (shared) -- (second);
    \foreach \column in {0.5,1.5,2.5,3.5}
        \foreach \row in {0.5,1.5,2.5,3.5}
            \node[locker] at (\column,\row) {};
    \node[customer] at (first) {};
    \node[customer] at (second) {};

    \node[below=3pt,inner sep=1pt,font=\footnotesize] at (0.5,0.5) {Depot};

    \draw[route] (-0.70,-0.42) -- (-0.38,-0.42);
    \node[anchor=west] at (-0.29,-0.42) {Vehicle route};
    \node[locker] at (1.94,-0.42) {};
    \node[anchor=west] at (2.08,-0.42) {Locker};
    \node[customer] at (3.48,-0.42) {};
    \node[anchor=west] at (3.62,-0.42) {Customer};
\end{tikzpicture}
\caption{An example of the grid route policy. Each square cell has diameter at most $2r$ and a locker at its center. The vehicle visits every locker, serves orders only there, and returns to the lower-left locker, its depot. For example, the two customers in the shaded cell share a locker, where their orders are served one at a time.}
\label{fig:fixed-route}
\end{wrapfigure}

\leavevmode When $r$ is of order $1-\rho$ or larger, all stable fixed-route policies with route length $L_r=O(1/r)$ attain the optimal order. If $r$ shrinks faster than $1-\rho$, every fixed-route policy fails to attain the optimal order $(1-\rho)^{-2}$, while a dynamic routing policy~\cite[Section~4.3]{BvR1991} attains it by choosing its route dynamically based on the addresses in each batch, with an appropriately chosen batch size.

The \emph{grid route policy} is one example of a fixed-route policy. The carrier chooses the service points and the route together on a grid of equal square cells of diameter at most $2r$. It uses the smallest even number of cells per side that meets this bound. A closed tour joins neighboring centers and visits each once. The vehicle serves each demand at the center of its cell, within distance $r$ of its address. In parcel delivery, these service points can be lockers or other pickup points, as illustrated in \Cref{fig:fixed-route}. In a mobile service program, they can be regional service centers where customers meet the vehicle.
\par
\endgroup

\subsection{Related Literature}

Our paper is most closely related to dynamic vehicle routing, which provides the model and mean system time benchmarks for our analysis. Out-of-home delivery and routing with neighborhoods study the choice of service locations and the associated travel savings. We also relate our work to flexibility design, which studies how to allocate limited flexibility, and to spatial queueing, which provides methods for analyzing system time along a prescribed route.

\paragraph{Dynamic vehicle routing.}
In dynamic routing, requests continue to arrive while a vehicle is on the road \cite{Psaraftis1980,Psaraftis1995,Pillac2013}. With no service location flexibility, heavy-traffic bounds relate mean system time to the travel required to serve waiting demands \cite{BvR1991,Bertsimas1993,Bertsimas1993multi,xu1994optimal}. Bertsimas and van Ryzin analyze tours through batches of a prescribed size~\cite[Section~4.3]{BvR1991}.
Bertsimas and van Ryzin~\cite[Section~4.4]{BvR1991} also analyze a policy proposed by Bartholdi and Platzman~\cite{BartholdiPlatzman1988} that uses a space-filling curve to order visits to outstanding demand locations. In our fixed-route policies, the vehicle repeats the same physical route on every trip and serves each demand at an assigned point within distance $r$ of its location.
The one-region divide-and-conquer policy of Pavone, Frazzoli, and Bullo takes a shortest tour through all demands waiting when the tour is planned. For a fixed service-time distribution with positive mean and finite second moment, they bound its mean system time by $O((1-\rho)^{-2})$ as $\rho\to1$~\cite[Theorem~4.2]{Pavone2011}.

\c{C}elik and Modiano~\cite{Celik2012} study a mobile collector that can receive a message within distance $r$ of its sender, with a fixed transmission time per message. For any fixed $r>0$, they show that their partitioning policy attains the optimal order $(1-\rho)^{-1}$ for mean delay. The collector cycles through the centers of square cells, serving each cell until it is empty. In our model, the service radius can also shrink as the load rises. Even then, the optimal mean system time is of smaller order than under exact-location service precisely when $r$ decreases more slowly than $1-\rho$. For $r$ of order $1-\rho$ or larger, every stable policy in our fixed-route class with an $r$-covering route of length $O(1/r)$ attains the optimal order.

\paragraph{Out-of-home delivery.}
Janinhoff et al.\ survey delivery systems in which customers can receive an order at a location other than their home \cite{Janinhoff2024}. Routing models let carriers choose among delivery options that customers rank \cite{Tilk2021,Dumez2021}. Using parcel delivery data, Krushinsky, Guo, and Claassen estimate the travel savings from alternative addresses and shared stops \cite{Krushinsky2021}. Carriers can influence customer choice through the sites and schedules of pickup trucks \cite{Glaeser2019} or through incentives to collect an order at a pickup point \cite{Galiullina2024}. The route savings from a locker network also depend on which customers use it \cite{Wang2025}. Ulmer and Streng allow delivery to a customer's preferred pickup station or another nearby station~\cite{UlmerStreng2019}. Our model allows any service point within the agreed distance and holds demand fixed to isolate the routing benefit of that choice. A given locker network would restrict this choice to its sites. In our construction, the carrier chooses the locker locations and the grid route together before demands arrive, so every customer has an assigned locker within distance $r$ (\Cref{fig:fixed-route}).

\paragraph{Traveling salesperson problem with neighborhoods.}
For a given batch, the goal is to find a shortest closed tour that intersects each demand's service disk \cite{Dumitrescu2003,Behdani2014,Bandyapadhyay2025}. Two feasible routes are a shortest tour through the demand locations and a covering route with parallel passes that bring the vehicle within service range of every point in the square \cite{ARKIN200025,Fekete}. For independent uniform locations in the unit square and shrinking radii, the shorter of these routes attains the optimal order of expected tour length \cite{CarlssonEtAl2016,Carlsson}. 

\paragraph{Flexibility design.}
In feature-based matching markets, a supplier can perform a task if its characteristics lie in the supplier's service region. The volume of this region measures the supplier's flexibility. When each supplier can serve at most one demand, allocating a fixed total service volume more evenly across suppliers can increase the expected number of demands fulfilled~\cite{Uniformity2026}. On ride-hailing and delivery platforms, a driver's pickup radius measures flexibility. Concentrating a small total radius budget on fewer drivers can fulfill more demands on average than an equal split~\cite{AmeenYu2026}. More broadly, work on flexibility design asks how to allocate flexibility across the two sides of a market \cite{Freund2024,Flexibility2026}. In our model the vehicle can serve many demands in succession. Moving one service point closer to the next can bring forward service for later arrivals as well.

\paragraph{Spatial queueing.}
Berman, Larson, and Chiu study where to base a mobile server on a network when a trip to one customer delays the response to the next \cite{Berman1985}. Fuhrmann and Cooper~\cite{FuhrmannCooper1985Continuum} study a server that collects demands arriving uniformly along a closed route and serves them upon returning to the origin. They show that mean system time exceeds the corresponding M/G/1 mean by one mean cycle duration. Under a fixed-route policy, the vehicle takes all demands waiting at departure, as in globally gated polling~\cite{BoxmaLevyYechiali1992}, and serves them along the route. Altman and Levy analyze stability and waiting times under a gated scan along a fixed path, with a separate round trip from the path to each customer~\cite[Section~5]{Altman1994}. Our fixed-route policies have the same gated-cycle structure, with service at assigned service points on the covering route and no local detours. Kavitha and Altman~\cite{KavithaAltman2012} compare routes for a wireless ferry that serves users in a two-dimensional region from assigned points along its path. A longer route can bring the ferry closer to users and shorten the time needed to exchange messages. In our model, a demand has the same on-site service requirement at any point within distance $r$ of its address. We characterize when a fixed-route policy attains the optimal order of mean system time among all admissible stable policies (\Cref{thm:fixed-route-scaling}).

\subsection{Notation and paper organization}
For $g>0$, we write $f=O(g)$ if $|f|/g$ stays bounded in the limit under consideration, and $f=\Theta(g)$ if $f/g$ stays between two positive constants. In heavy traffic, $r\gg1-\rho$ means that $r/(1-\rho)\to\infty$ as $\rho\to1$.
\Cref{sec:model} introduces the model and states the main results. \Cref{sec:analysis} proves them, and \Cref{sec:conclusion} discusses how pickup distances might vary across customers.

\section{Main Results}\label{sec:model}

\subsection[Model setup]{Model setup}
Demands arrive according to a Poisson process of rate $\lambda>0$. Each location is drawn uniformly from the unit square $[0,1]^2$, independently of the other locations and the arrival times. A single vehicle travels at speed $v=1$. On-site service times are independent copies of a nonnegative random variable $S$, independent of arrivals and locations, with mean $\bar s>0$ and finite second moment. The vehicle serves one demand at a time without moving, and each on-site service runs uninterrupted to completion. Overlapping neighborhoods permit several demands to be served sequentially at the same location.

Write $\rho=\lambda\bar s$ for the on-site service load. We require $\rho<1$, a necessary condition for stability. The arrival rate and the service-time distribution are fixed during operation. When we compare different loads, we may vary the arrival rate $\lambda$ and mean service time $\bar s$ together or hold either one fixed. At all loads, $\lambda$ stays bounded above and away from zero and $\E{S^2}\le K\bar s^2$ for the same finite constant $K>0$. Under these bounds, $\bar s=\rho/\lambda$ stays bounded above and away from zero as $\rho\to1$. This includes the classical limit in dynamic routing and spatial polling in which the service-time distribution is fixed and $\lambda\to1/\bar s$~\cite{Pavone2011,Kroese1997}.

Fix a service radius $r>0$ before operation, where $r$ may depend on the system load $\rho$, and keep it fixed over time. Every demand and every policy in the comparison uses this same radius. A demand at $\xbf\in[0,1]^2$ can be served at any point in
\[
 \{\ybf\in[0,1]^2:\norm{\ybf-\xbf}_2\le r\}.
\]

When it takes a new assignment, the vehicle selects one or more waiting demands and fixes a route from its current location to serve them. Once it completes the route and the assigned services, it may take a new assignment where the route ends, wait, or relocate within the square. An \emph{admissible policy} makes assignment and routing decisions using only information observed so far and independent random choices, without information about future arrivals or the service durations of demands whose service has not begun.

Let $N(t)$ be the number of outstanding demands at time $t$, including any demand in service. A policy is stable if arrivals, vehicle motion, and service admit a joint stationary regime in which $N(t)<\infty$ almost surely and every arrival eventually completes service.

A demand's system time is the duration from arrival to service completion. Let $\pi$ specify an admissible stable policy and one of its stationary regimes. Write $T_\pi$ for the mean system time of a typical arrival in that regime. The optimal mean system time $T_*$ is the infimum of $T_\pi$ over all admissible stable policies and their stationary regimes at the same arrival rate, service-time distribution, and radius.

\subsection[Optimal mean system time in heavy traffic]{Optimal mean system time in heavy traffic}\label{sec:mainresults}
We first characterize how the optimal mean system time $T_*$, the infimum of $T_\pi$ over all admissible stable policies, scales with the service radius $r$ in the following theorem.

\begin{theorem}[Heavy-traffic scaling of optimal mean system time]\label{thm:scaling}
For all $r>0$ and $c\le\rho<1$, where $c\in(0,1)$ is any fixed constant, the optimal mean system time satisfies
\[
 T_*=\Theta\!\left(\frac{1+\min\{1/r,\,1/(1-\rho)\}}{1-\rho}\right).
\]
\end{theorem}

To compare service radii that shrink at different rates as the load rises, we consider the family
\begin{equation}\label{eq:service-radius}
 r=(1-\rho)^\alpha, \qquad \alpha\ge0.
\end{equation}
At $\alpha=1$, the radius equals $1-\rho$, the threshold scale in \Cref{thm:scaling}. As the load rises, the radius shrinks more slowly for $0<\alpha<1$ and faster for $\alpha>1$. At any given load, $\alpha=0$ gives $r=1$, and $r\to0$ as $\alpha\to\infty$.
Let $\Pi(\alpha)$ denote all admissible stable policies paired with their stationary regimes, and write
\[
 T_*(\alpha)=\inf_{\pi\in\Pi(\alpha)} T_\pi(\alpha).
\] 
When an optimal policy exists, $T_*(\alpha)$ is the average time its demands spend waiting and receiving service. Substituting $r=(1-\rho)^\alpha$ into \Cref{thm:scaling} gives the following corollary.

\begin{corollary}[Optimal scaling for power-law radii]\label{cor:power-radius}
For any fixed $\alpha\ge0$, let $r=(1-\rho)^\alpha$. As $\rho\to1$,
\begin{equation}
 T_*(\alpha)=
 \begin{cases}
  \Theta\bigl((1-\rho)^{-(1+\alpha)}\bigr),&0\le\alpha\le1,\\
  \Theta\bigl((1-\rho)^{-2}\bigr),&\alpha>1.
 \end{cases}
 \label{eq:alphascaling}
\end{equation}
\end{corollary}

For comparison, consider a single server working in arrival order without travel. By the Pollaczek--Khinchine formula for this M/G/1 queue~\cite[Section~3.5]{BertsekasGallager1992}, its mean system time is
\begin{equation}\label{eq:queue-benchmark}
 T_{\rm queue}=\bar s+\frac{\lambda\E{S^2}}{2(1-\rho)}.
\end{equation}

\begin{remark}[Single-stop and exact-location service]
The center of the unit square is within distance $r$ of every location when $r\ge1/\sqrt{2}$, half the diagonal. Serving demands there in arrival order attains the M/G/1 benchmark in \eqref{eq:queue-benchmark}, so $T_*=T_{\rm queue}=\Theta((1-\rho)^{-1})$ as $\rho\to1$~\cite[Section~3.5]{BertsekasGallager1992}. At $r=0$, the vehicle must visit every demand location, and $T_*=\Theta((1-\rho)^{-2})$ as in the classical dynamic traveling repairman problem~\cite{BvR1991}.
\end{remark}

\subsection[Fixed-route policies]{Fixed-route policies}\label{sec:fixed-route-policies}
In this subsection, we introduce the class of fixed-route policies ($\FR$). We characterize when the service radius is large enough for a policy in $\FR$ to attain the optimal heavy-traffic order of mean system time and when every fixed-route policy falls short.

An $r$-covering route is a continuous closed curve of finite length in the unit square that comes within distance $r$ of every point in the square. Before operation, a policy in $\FR$ chooses such a route and designates a point on it as the depot. A departure is the start of a cycle at the depot. The vehicle takes all demands waiting at departure and chooses its direction by any admissible rule, for example always traveling forward or tossing an independent fair coin.

\begin{figure}[h]
\centering
\begin{tikzpicture}[
    x=1.10cm,y=1.10cm,font=\small,line join=round,line cap=round,
    traveled/.style={draw=darkgreen,line width=1.4pt},
    remaining/.style={draw=black!50,dashed,line width=0.75pt},
    direction/.style={traveled,-{Stealth[length=7pt,width=6pt]}},
    locker/.style={rectangle,draw=black!55,fill=white,line width=0.7pt,
                   minimum size=3.8pt,inner sep=0pt}
]
\foreach \panel/\offset in {1/0,2/5.5}{
    \begin{scope}[shift={(\offset,0)}]
        \coordinate (depot) at (0.5,0.5);
        \coordinate (demand) at (1.14,2.78);
        \coordinate (service) at (1.5,2.5);
        \def\forwardarc{(depot) -- (0.5,3.5) -- (1.5,3.5) -- (service)}
        \def\reversearc{(depot) -- (3.5,0.5) -- (3.5,3.5) -- (2.5,3.5) --
            (2.5,1.5) -- (1.5,1.5) -- (service)}
        \fill[darkgreen!12] (1,2) rectangle (2,3);
        \draw[xstep=1.10cm,ystep=1.10cm,black!25,line width=0.35pt]
            (0,0) grid (4,4);
        \draw[black!40,line width=0.45pt] (0,0) rectangle (4,4);
        \ifnum\panel=1
            \draw[remaining] \reversearc;
            \draw[traveled] \forwardarc;
            \draw[direction] (0.5,1.68) -- (0.5,2.18);
            \node[above=3pt] at (2,4) {Forward distance $t(x)$};
        \else
            \draw[remaining] \forwardarc;
            \draw[traveled] \reversearc;
            \draw[direction] (3.5,1.68) -- (3.5,2.18);
            \draw[direction] (0.68,0.5) -- (1.18,0.5);
            \node[above=3pt] at (2,4) {Reverse distance $L_r-t(x)$};
        \fi
        \foreach \column in {0.5,1.5,2.5,3.5}
            \foreach \row in {0.5,1.5,2.5,3.5}
                \node[locker] at (\column,\row) {};
        \draw[black!65,densely dotted,line width=0.7pt] (demand) -- (service);
        \fill[black] (demand) circle[radius=1.5pt] node[above left=1pt] {$x$};
        \node[locker,draw=darkgreen,line width=1pt,minimum size=5.5pt]
            at (service) {};
        \node[right=2pt] at (service) {$p$};
        \node[rectangle,fill=black,inner sep=0pt,minimum size=5.5pt] at (depot) {};
        \node[below=1pt,inner sep=1pt,font=\footnotesize] at (depot) {Depot};
    \end{scope}
}
\end{tikzpicture}
\caption{An example of a fixed-route policy: the grid route policy with $r=\sqrt{2}/8$ (see \Cref{rem:grid-route} for details). The solid paths from the depot to the same assigned locker $p=\Gamma(t(x))$ have lengths $t(x)$ and $L_r-t(x)$.}
\label{fig:route-directions}
\end{figure}
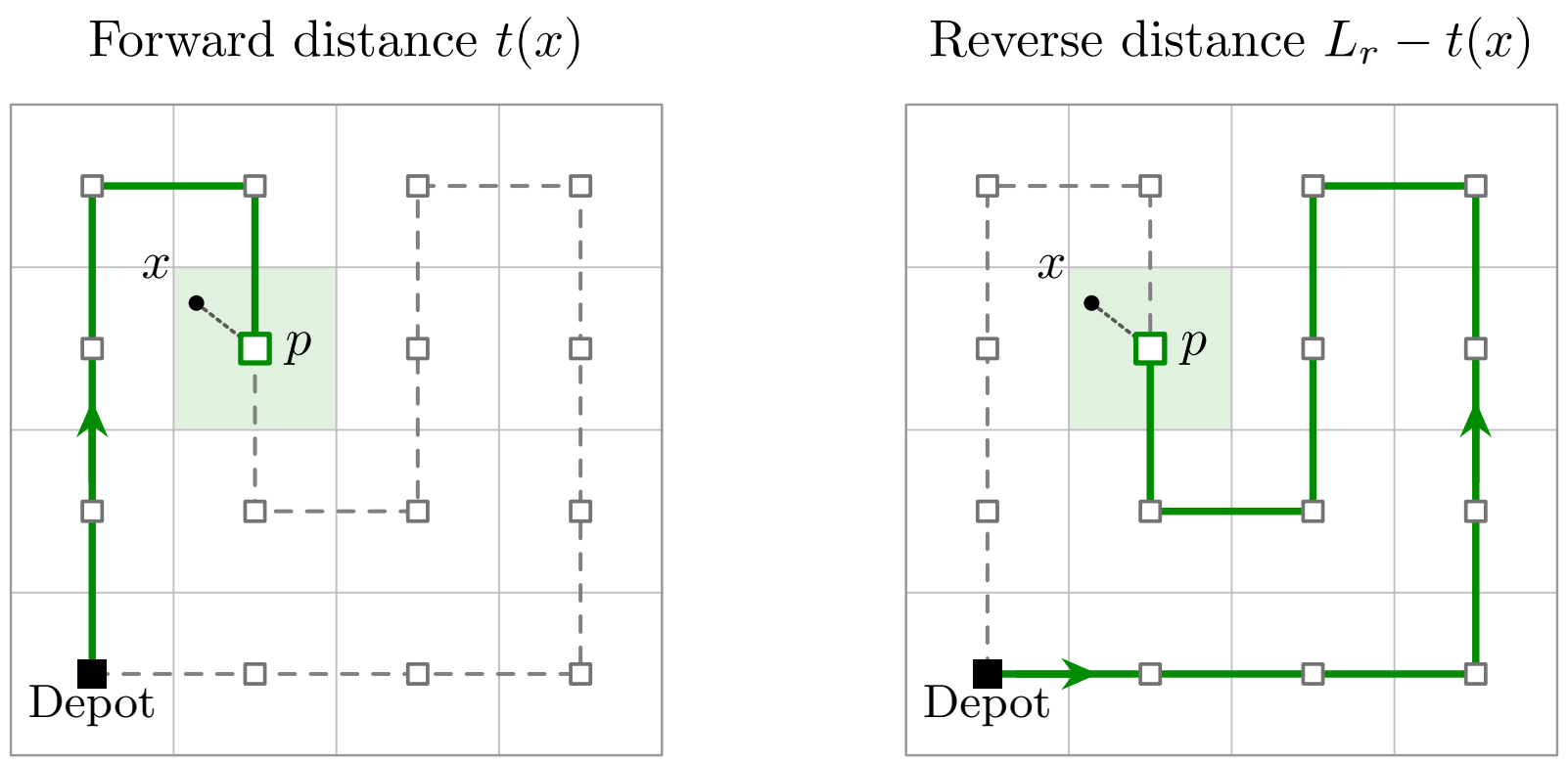

Write $L_r>0$ for the length of the chosen route. Let $\Gamma:[0,L_r]\to[0,1]^2$ trace it by arclength from the depot in a fixed reference direction, with $\Gamma(0)=\Gamma(L_r)$. Before operation, fix a measurable function $t:[0,1]^2\to[0,L_r)$ that assigns each location $x$ a route parameter $t(x)$ satisfying $\norm{\Gamma(t(x))-x}_2\le r$. 
The vehicle travels at unit speed and stops to serve $x$ after traveling $t(x)$ forward or $L_r-t(x)$ in reverse, as shown in \Cref{fig:route-directions}. Thus a demand assigned $t(x)=0$ is served at the depot before travel begins in the forward direction or after the vehicle returns in reverse. If several demands are assigned to the same visit, the vehicle serves them in the order in which they arrived. The cycle ends at the depot after the full traversal and all assigned services. The next cycle starts immediately, even with an empty batch.

When the vehicle always travels forward, it serves all demands waiting at departure in the order of their assigned positions $t(x)$ along the route, as in globally gated polling. Under our service and travel assumptions, the waiting-time formula of Boxma, Levy, and Yechiali~\cite[Eq.~(2.14)]{BoxmaLevyYechiali1992} yields the forward mean system time for finitely many service points. For continuously distributed service positions, the workload formula of Kavitha and Altman~\cite[Eq.~(5)]{KavithaAltman2012} yields the same expression. The following theorem gives the exact mean system time for either fixed direction in terms of the route length and the mean forward distance from the depot to a demand's assigned service point. We use this identity in \Cref{thm:fixed-route-scaling} to compare fixed-route policies with the optimum.

\begin{theorem}[Exact mean system time for fixed-route policies]\label{thm:fixedradius}
Fix $r>0$ and an $r$-covering route of length $L_r>0$. Let $X$ be a uniformly distributed address in the unit square. For $0<\rho<1$, let $T_{\mathrm{FR}}$ denote the stationary mean system time.
\begin{equation}\label{eq:fixed-direction-T}
 T_{\mathrm{FR}}=T_{\rm queue}+\frac{L_r}{1-\rho}
 +\begin{cases}
    \E{t(X)}-L_r/2, & \text{if the vehicle always travels forward},\\
    L_r/2-\E{t(X)}, & \text{if the vehicle always travels in reverse}.
  \end{cases}
\end{equation}
\end{theorem}

\begin{remark}[Mean cycle duration and the fair-coin rule]
The first term, $T_{\rm queue}$, is the mean system time in the non-spatial M/G/1 queue in \eqref{eq:queue-benchmark}. Each cycle requires $L_r$ units of travel, and the vehicle spends a fraction $1-\rho$ of its time traveling. Thus the second term, $L_r/(1-\rho)$, is the mean time between departures, including on-site service. For example, a stop a quarter of a lap from the depot going forward is three quarters of a lap away in reverse. If the travel distance to a typical demand's service point averages less than half a lap going forward, mean system time is shorter in that direction. With an independent fair coin at each departure, the two distances to each service point are equally likely and average to $L_r/2$. The third term therefore vanishes and
\begin{equation}
 T_{\mathrm{FR}}=T_{\rm queue}+\frac{L_r}{1-\rho}.
 \label{eq:fixedT}
\end{equation}
Note that both \eqref{eq:fixed-direction-T} and \eqref{eq:fixedT} are exact for every $0<\rho<1$, including intermediate loads. Under the fair-coin rule, \eqref{eq:fixedT} explicitly characterizes how mean system time depends on route length $L_r$.
\end{remark}

The following theorem compares the stationary mean system time under any admissible tie-breaking rule for the direction at departure\footnote{The direction may depend on the waiting demands and past cycles, using independent randomness if desired, but not on future arrivals or the service durations of demands whose service has not begun.} with the optimal order in \Cref{thm:scaling}.

\begin{theorem}[Order optimality of fixed-route policies]\label{thm:fixed-route-scaling}
For any fixed $\alpha\ge0$ and any stable policy in $\FR$ whose prescribed route at radius $r=(1-\rho)^\alpha$ has length $L_r=O(1/r)$ as $\rho\to1$,
\[
 T_{\text{FR}}(\alpha)=\Theta\!\left((1-\rho)^{-(1+\alpha)}\right),
\]
which matches the optimal order in \eqref{eq:alphascaling} if and only if $0\le\alpha\le1$.
When $r/(1-\rho)\to0$ as $\rho\to1$, the ratio $T_{\mathrm{FR}}/T_*$ tends to infinity for every stable policy in $\FR$, regardless of its route length.
\end{theorem}

\begin{figure}[h]
\centering
\definecolor{optline}{HTML}{1F5FBF}
\definecolor{fixline}{HTML}{D96C00}
\definecolor{gainregion}{HTML}{E4F5EF}
\definecolor{nogainregion}{HTML}{F0EFEC}
\begin{tikzpicture}[plot annotation/.style={font=\fontsize{9.5}{11}\selectfont}]
\begin{axis}[
    width=0.74\textwidth,
    height=0.46\textwidth,
    xmin=0, xmax=2.2,
    ymin=0.7, ymax=3.5,
    axis lines=left,
    axis on top,
    axis line style={-{Stealth[length=4pt]}, black!60, line width=0.4pt},
    tick style={black!60},
    xtick={0,0.5,1,1.5,2},
    ytick={1,1.5,2,2.5,3},
    xlabel={$\alpha$, where the radius $r=(1-\rho)^\alpha$},
    ylabel={Exponent $\kappa$ for mean system time\\$T=\Theta\bigl((1-\rho)^{-\kappa}\bigr)$},
    ylabel style={align=center},
    label style={font=\small},
    tick label style={font=\small},
    clip=false,
]
\fill[gainregion] (axis cs:0,0.7) rectangle (axis cs:1,3.5);
\fill[nogainregion] (axis cs:1,0.7) rectangle (axis cs:2.2,3.5);
\node[plot annotation, align=center] at (axis cs:0.5,3.24)
    {Flexibility succeeds\\Optimal exponent $<2$};
\node[plot annotation, align=center] at (axis cs:1.55,3.24)
    {Flexibility fails\\Optimal exponent $=2$};
\draw[black!65, dashed, line width=0.7pt]
    (axis cs:1,0.7) -- (axis cs:1,3.5);
\addplot[black!50, densely dotted, line width=1pt]
    coordinates {(0,2) (1,2)};
\node[plot annotation, align=center, anchor=south]
    (benchmarklabel) at (axis cs:0.46,2.30) {Benchmark with\\no flexibility};
\draw[black!65, line width=0.7pt, -{Stealth[length=4pt]}]
    (benchmarklabel.south) -- (axis cs:0.32,2);
\addplot[fixline, line width=3.6pt, dash pattern=on 6pt off 4pt]
    coordinates {(0,1) (2.2,3.2)};
\node[plot annotation, align=center, anchor=north]
    (fixedlabel) at (axis cs:1.91,2.56)
    {Fixed-route policy\\$T_{\mathrm{FR}}$, $\kappa=1+\alpha$};
\draw[fixline, line width=0.7pt, -{Stealth[length=4pt]}]
    (axis cs:1.85,2.85) -- (fixedlabel.north);
\addplot[optline, line width=1.3pt]
    coordinates {(0,1) (1,2) (2.2,2)};
\node[plot annotation, align=center, anchor=north]
    (optimallabel) at (axis cs:1.62,1.65)
    {Optimum $T_*(\alpha)$\\$\kappa=1+\min\{\alpha,1\}$};
\draw[optline, line width=0.7pt, -{Stealth[length=4pt]}]
    (axis cs:1.52,2) -- (optimallabel.north);
\fill[optline] (axis cs:0,1) circle[radius=2.2pt];
\fill[optline] (axis cs:1,2) circle[radius=2.2pt];
\node[plot annotation, anchor=north west] at (axis cs:0.04,0.98)
    {M/G/1 optimum ($\alpha=0$)};
\node[font=\small, anchor=south, inner sep=2pt] at (rel axis cs:0,1) {$\kappa$};
\node[font=\small, anchor=west, inner sep=2pt] at (rel axis cs:1,0) {$\alpha$};
\end{axis}
\end{tikzpicture}
\caption{Heavy-traffic exponent $\kappa$ of mean system time for $r=(1-\rho)^\alpha$ as $\rho\to1$. The solid blue curve shows $\kappa$ for the optimal mean system time $T_*(\alpha)$ (\Cref{cor:power-radius}). The dashed orange curve shows the exponent $1+\alpha$ for any stable policy in $\FR$ with $L_r=O(1/r)$ (\Cref{thm:fixed-route-scaling}). The two exponents agree for $0\le\alpha\le1$. The dotted line gives the benchmark with no service location flexibility ($r=0$)~\cite{BvR1991}.}
\label{fig:alpha-exponent}
\end{figure}
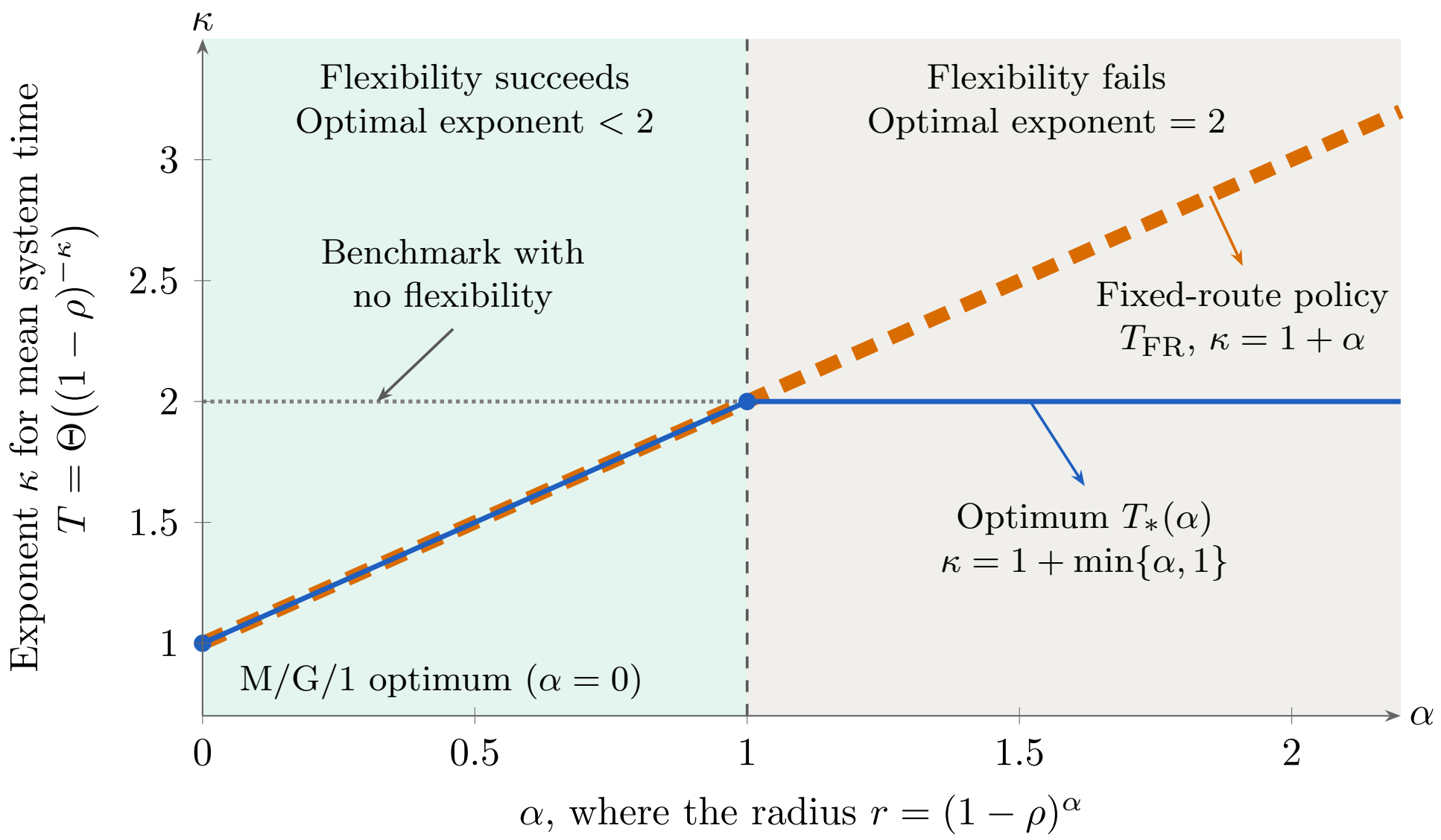

\Cref{fig:alpha-exponent} compares the optimal exponent with that of fixed-route policies satisfying $L_r=O(1/r)$.

\FloatBarrier

\begin{remark}[The grid route and fixed pickup points]\label{rem:grid-route}
The grid route in \Cref{fig:fixed-route,fig:route-directions} extends to any service radius $r>0$. Let $k$ be the smallest positive even integer with $1/k\le\sqrt{2}\,r$. Divide the unit square into $k^2$ equal square cells. If an address lies on a shared edge or corner, choose one of the adjoining cells uniformly at random before operation and assign the address to that cell.\footnote{A demand falls on a grid boundary with probability zero, so ties can be broken arbitrarily.} The vehicle serves each demand at the center of its assigned cell. The distance from the center to a corner is $1/(\sqrt{2}k)\le r$, so the center can serve the entire cell. In the forward direction, the vehicle leaves the depot at the lower-left cell center and visits the centers above the bottom row in alternating passes up and down the columns. It then returns to the depot along the bottom row.

In \Cref{fig:route-directions}, $r=\sqrt{2}/8$ and $k=4$. The route visits all $16$ lockers and has length $L_r=4$. The demand at $x$ is served at the same locker $p$, the center of its cell, in both directions. In general, the grid route has length $L_r=k=O(1/r)$ as $r\to0$. Thus, for $r=(1-\rho)^\alpha$, the grid route policy attains the optimal order of mean system time for $0\le\alpha\le1$ by \Cref{thm:fixed-route-scaling}.
\end{remark}

\begin{remark}[Intuition for the radius threshold]
In steady state, on-site service occupies a fraction $\rho$ of the vehicle's time, leaving at most $(1-\rho)/\lambda$ units of travel per completed demand on average. This allowance has order $1-\rho$ under our bounds on $\lambda$. A shortest tour through $q$ independent uniform addresses has expected length of order $\sqrt q$, or order $1/\sqrt q$ per demand~\cite{Beardwood_Halton_Hammersley_1959}. A batch of order $(1-\rho)^{-2}$ therefore brings the travel per demand down to order $1-\rho$, while the expected total tour length has order $(1-\rho)^{-1}$. A fixed-route policy can use an $r$-covering route of length of order $1/r$ as $r\to0$, such as the grid route. This length does not depend on the batch size. When $r\gg1-\rho$, this is of smaller order than the expected length of that tour, suggesting that fewer demands need to share its travel.
\end{remark}

\section{Analysis}\label{sec:analysis}

We now provide the theoretical backbone of our main results in \Cref{sec:model}. The section is organized as follows. In \Cref{sec:gates}, we provide two propositions required to prove the theorems. In \Cref{prop:gatedmoments}, we relate the mean system time of each policy in $\FR$ that has a stationary regime to its batch size at departure. We also use the travel needed to serve recent arrivals to derive a lower bound for every admissible stable policy, in \Cref{thm:unrestricted}. Using \Cref{prop:gatedmoments}, we first prove \Cref{thm:fixedradius} in \Cref{sec:fixedroute} to obtain one of the upper bounds required for \Cref{thm:scaling}. Furthermore, we prove \Cref{thm:scaling} in \Cref{sec:upper} by obtaining the other upper bound from tours through the addresses in each batch constructed in \cite{BvR1991} and showing that the smaller of the two matches the lower bound in order. Finally, we compare the fixed-route and optimal orders in \Cref{sec:proof-fixed-route-scaling} to prove \Cref{thm:fixed-route-scaling}.

\subsection{Fixed-route policies and a lower bound for all policies}\label{sec:gates}\label{FixedShapeDC}\label{sec:fixed-radius-lower}\label{sec:lower}\label{sec:temporal-cohort-lower-bound}

For any policy in the class $\FR$ defined in \Cref{sec:fixed-route-policies}, a cycle $i$ denotes the time from departure $i$ to departure $i+1$ at the depot, with departures indexed by $i\in\mathbb{N}$.

At departure $i$, the number of available demands is denoted by $Q_i$. At the time of departure, the vehicle takes all $Q_i$ waiting demands and fixes their service order according to its direction rule. For $i\ge1$, batch $i$ consists of the demands that arrived during cycle $i-1$. By the definition of $\FR$ class policies, travel takes $L_r$ time units in every cycle, even when $Q_i=0$. 
Let $S_{ij}$ be the service time of the $j$th demand served in cycle $i$. Then cycle $i$ has duration
\[
\tau_i=L_r+\sum_{j=1}^{Q_i}S_{ij} \,.
\]

Each departure time is a stopping time, so arrivals following a departure form a Poisson process of rate $\lambda$ independent of the history up to that departure. Given the route and service times, $Q_{i+1}$ counts all arrivals in an interval of known length $\tau_i$. Therefore, its conditional distribution depends only on $\tau_i$, so
\begin{equation}
 Q_{i+1}\mid \tau_i \sim \operatorname{Pois}(\lambda \tau_i).
 \label{eq:gatedrecursion}
\end{equation}

The following proposition establishes the relationship between the stationary batch-size moments under any direction rule and mean system time.

\begin{proposition}[Stability and mean system time]\label{prop:stability}\label{prop:gatedmoments}
For every $0<\rho<1$ and every policy in $\FR$, the batch sizes $(Q_i)_{i\ge1}$ form a positive recurrent Markov chain. Let $Q$ denote the batch size under its stationary distribution. Then $\E{Q^2}<\infty$, and the ratio $R_Q=\E{Q^2}/\E{Q}$ satisfies

\begin{equation}
 R_Q=\frac{\lambda^2\E{S^2}+(1+\rho)\lambda L_r}{1-\rho^2}+1.
 \label{eq:gatesecond}
\end{equation}
The fixed direction and fair coin rules are stable. In any stationary regime of a policy in $\FR$, let $h_j$ be the travel time from departure to the $j$th service, excluding all service durations, for $1\le j\le Q$. The stationary mean system time is finite and satisfies
\begin{equation}
 T_{\rm FR}=\frac{R_Q-1}{2\lambda}+\frac{\bar s(R_Q+1)}{2}
   +\frac{\E{\sum_{j=1}^Q h_j}}{\E{Q}},
 \label{eq:rewardbound}
\end{equation}
where expectations on the right are taken at stationary departures.
\end{proposition}
 The relationship between $T_\FR$ and the stationary batch sizes becomes clearer when the system time result of \Cref{prop:gatedmoments} is interpreted in the following way. Consider a typical demand under the stationary distribution.
Each batch of size $q$ serves $q$ demands, so a typical demand belongs to such a batch with probability $q\Pp(Q=q)/\E{Q}$.  Therefore, the expected batch size of a typical demand is given by $\sum_{q\ge1}q^2\Pp(Q=q)/\E{Q}=R_Q$.

Finally, the following proposition gives a lower bound on mean system time for every admissible stable policy.

\begin{proposition}[Uniform lower bound]\label{thm:fixed-radius-lower}\label{thm:unrestricted}
For every $r>0$ and $0<\rho<1$, every admissible stable policy has stationary mean system time $T\ge T_{\rm queue}$.  There exist constants $r_0,c_0>0$,  independent of problem parameters, such that for $0<r\leq r_0$,\footnote{A single service point covers an area of at most $\pi r^2$. When $r$ is small, the vehicle must therefore travel to serve demands spread across the square. For each $\alpha>0$, the prescribed radius $r=(1-\rho)^\alpha$ eventually lies below $r_0$ as $\rho\to1$.} every policy also satisfies
\begin{equation}
 T\ge c_0\min\left\{\frac{1}{r(1-\rho)},\frac{\lambda}{(1-\rho)^2}\right\}.
 \label{eq:rate-explicit-lower}
\end{equation}
\end{proposition}

\ifhideproofs\else
\subsubsection[Proof of Proposition~\ref*{prop:stability}]{Proof of Proposition~\ref*{prop:stability}}

 We first show that $(Q_i)_{i\geq1}$ admits a stationary distribution. Given $Q_i=q$, the cycle time $\tau_i$ consists of $L_r$ and $q$ independent service times, regardless of the demand locations and the chosen direction. By \eqref{eq:gatedrecursion}, $(Q_i)_{i\geq 1}$ satisfies the Markov property as the distribution of the next batch size, $Q_{i+1}$, depends on the past only through $q$. Since every policy in $\FR$ has $L_r>0$, the Poisson distribution gives positive probability to every next batch size. Thus $(Q_i)_{i\ge1}$ is an irreducible, aperiodic Markov chain.

Given $Q_i=q$, the conditional mean and second factorial moment of $Q_{i+1}$ are
\begin{equation}
 \E{Q_{i+1}\mid Q_i=q}=\lambda L_r+\rho q,
 \label{eq:gate-conditional-mean}
\end{equation}
\begin{equation}
 \E{Q_{i+1}(Q_{i+1}-1)\mid Q_i=q}
 =(\lambda L_r+\rho q)^2+\lambda^2q\Var(S),
 \label{eq:gate-conditional-factorial}
\end{equation}
where the second equality uses the Poisson factorial moment $(\lambda\tau_i)^2$: conditional on $Q_i=q$, $\lambda\tau_i$ has mean $\lambda L_r+\rho q$ and variance $\lambda^2q\Var(S)$. 
Adding these two identities and subtracting $q^2$ gives
\[
 \E{Q_{i+1}^2-Q_i^2\mid Q_i=q}
 =-(1-\rho^2)q^2+O(q+1).
\]
For fixed system parameters, this is at most $-(1-\rho^2)q^2/2$ for all sufficiently large $q$, since $\rho<1$. For the finitely many smaller batch sizes, the conditional second moments are bounded. From each of these batch sizes, the next batch is empty with positive probability. The Foster--Lyapunov criterion applied to the Lyapunov function $V(q)=1+q^2$ gives positive recurrence and $\E{Q^2}<\infty$~\cite[Theorem~14.0.1]{MeynTweedie1993}.

 We will now prove \eqref{eq:gatesecond} under the stationary distribution of batch sizes. \eqref{eq:gate-conditional-mean} implies that
\begin{equation*}
    \E{Q_{i+1}} = \lambda L_r + \rho\E{Q_i}.
\end{equation*}
Under stationarity, we have that $\E{Q_{i+1}}=\E{Q_i}$, so we obtain $(1-\rho)\E{Q}=\lambda L_r$.  Applying the same argument to \eqref{eq:gate-conditional-factorial} and the stationarity condition gives
\[
 (1-\rho^2)\E{Q(Q-1)}
 =\bigl(\lambda^2\E{S^2}+2\rho\lambda L_r\bigr)\E{Q}+(\lambda L_r)^2,
\]
where $\rho^2+\lambda^2\Var(S)=\lambda^2\E{S^2}$. Dividing by $\E{Q}$ and using $(\lambda L_r)^2/\E{Q}=(1-\rho)\lambda L_r$ proves \eqref{eq:gatesecond}.

We now establish stationarity for the fixed direction and fair coin rules. If the system starts at a departure with the stationary batch size and independent and uniformly distributed locations, the next batch size has the same stationary distribution, and conditional on its size the new arrivals again have independent uniform locations. Since the direction is fixed or chosen by an independent fair coin, the sequence of cycles is also stationary. The mean cycle duration is positive and finite, because under the stationary distribution at departures, the mean cycle duration satisfies $0<\E{\tau_i}=L_r+\bar s\E{Q}<\infty$. To observe the vehicle at a typical time, we sample a cycle in proportion to its duration, weighting this distribution by $\tau_i/\E{\tau_i}$ (see~\cite[Chapter~$8$, Eq.~$(4.2^\circ)$]{Thorisson2000}). We then choose time zero uniformly within that cycle. By \cite[Chapter~$8$, Theorem~4.1]{Thorisson2000}, it follows that the process viewed from this time is stationary.

Now take any stationary regime of a policy in $\FR$. Its batch size at departure has the unique stationary distribution established above. We split each demand's system time at the start of the cycle that serves it. Each demand contributes once to each part, so its mean contribution is the expected total per cycle divided by $\E{Q}$. A demand arriving $u$ time units into cycle $i$ waits $\tau_i-u$ for the next departure. The cycle duration is independent of the Poisson arrival process after departure, so the expected total wait of arrivals during that cycle, conditional on $\tau_i$, is
\begin{equation}
    \lambda\int_0^{\tau_i}(\tau_i-u)\,\diff u=\frac{\lambda\tau_i^2}{2}.
\end{equation}
Dividing its expectation by the mean number $\E{Q}=\lambda\E{\tau_i}$ of arrivals per cycle gives the mean wait per demand,
\[
 \frac{\lambda\E{\tau_i^2}/2}{\E{Q}}
 =\frac{\E{\tau_i^2}}{2\E{\tau_i}}
 =\frac{R_Q-1}{2\lambda},
\]
where the second equality uses $\lambda^2\E{\tau_i^2}=\E{Q(Q-1)}$, obtained by averaging the Poisson factorial moment in \eqref{eq:gatedrecursion} at stationary departures.

After departure, the $j$th demand experiences $h_j$ units of travel and $j$ services, including its own. The service order is fixed before these durations are observed, so summing over the batch and dividing the expected total by $\E{Q}$ gives the mean time per demand after departure,
\[
 \frac{\E{\sum_{j=1}^Q(h_j+j\bar s)}}{\E{Q}}
 =\frac{\E{\sum_{j=1}^Q h_j}}{\E{Q}}+\frac{\bar s(R_Q+1)}{2},
\]
where the equality uses $\sum_{j=1}^Q j=Q(Q+1)/2$. Adding the two contributions proves \eqref{eq:rewardbound}. Each term is finite because $R_Q<\infty$ and $0\le h_j\le L_r$.\qed

\subsubsection[Proof of Proposition~\ref*{thm:fixed-radius-lower}]{Proof of Proposition~\ref*{thm:fixed-radius-lower}}

Fix the radius $r>0$ and a stable policy in a stationary regime. For this lower bound, we allow the vehicle to revise its route after starting an assignment to serve new arrivals, under the same information and service assumptions. It suffices to consider policies with finite mean system time $T$. Write $d=\E{N(0)}=\lambda T$ by Little's law.

We first show that $T\ge T_{\rm queue}$, following the workload comparison in \cite[proof of Theorem~3.1]{Altman1994}. Let $W$ be the total on-site service remaining at time zero. In stationarity, services start at rate $\lambda$ because every arrival eventually enters service. The mean number in service is therefore $\lambda\bar s=\rho$. This leaves $d-\rho$ demands waiting on average. Whether a demand is still waiting is independent of its own unobserved service time, so waiting demands contribute $\bar s(d-\rho)$ to $\E{W}$. A service of length $s$ contributes $s^2/2$ to the time integral of its remaining duration. Thus the mean remaining duration of the service in progress, counting zero when the vehicle is idle or traveling, is $\lambda\E{S^2}/2$. Adding the two contributions gives
\begin{equation*}
\E{W}=\bar s(d-\rho)+\frac{\lambda\E{S^2}}{2}.
\end{equation*}

Give the vehicle and the stationary queue in \eqref{eq:queue-benchmark} the same arrivals and service requirements. The queue's remaining work at time zero is at most $W$. Indeed, a busy queue was empty just before its current busy period and has worked continuously since, while the vehicle could have supplied at most as much service. Applying the preceding identity to the queue with $d=\lambda T_{\rm queue}$ gives $\lambda\bar s(T-T_{\rm queue})\ge0$.

Since the vehicle stays still during service and travels at unit speed, its path $\Gamma$ over any deterministic interval of length $h$ satisfies
\begin{equation}
 \E{|\Gamma|}\le(1-\rho)h.
 \label{eq:fixed-radius-travel-budget}
\end{equation}

To see how much travel recent arrivals require, divide the square into cells. Each cell that receives arrivals but has none of them served contributes at least one demand to the backlog.

\begin{lemma}[Area reached while serving recent arrivals]\label{lem:cohort-area}
Suppose $d=\lambda T\ge1$. Let $\Gamma$ be the vehicle's path during $[-64T,0]$ and set $a=r+1/(4\sqrt d)$. Then
\begin{equation}
 \E{\operatorname{area}\{x\in[0,1]^2:\operatorname{dist}(x,\Gamma)\le a\}}
 \ge \frac12.
 \label{eq:cohort-area}
\end{equation}
\end{lemma}

For any rectifiable path $\Gamma$ of length $|\Gamma|$ and any $a\ge0$, the tube-area bound in \cite[Eq.~(4.9)]{Steele1990} gives
\begin{equation}
 \operatorname{area}\{x\in[0,1]^2:\operatorname{dist}(x,\Gamma)\le a\}
 \le 2a|\Gamma|+\pi a^2.
 \label{eq:path-area}
\end{equation}

Take $r_0=1/16$ and suppose $r\le r_0$. First consider $d\ge1$. The radius $a$ in \Cref{lem:cohort-area} is at most $5/16$, so $\pi a^2<1/3$. Comparing the two area estimates gives
\[
 \frac12\le 2a\,\E{|\Gamma|}+\pi a^2
          \le 128a(1-\rho)T+\frac13,
\]
where the first inequality follows from \eqref{eq:cohort-area} and \eqref{eq:path-area}, and the second uses \eqref{eq:fixed-radius-travel-budget} over $64T$ time units. Since $aT=rT+\tfrac14\sqrt{T/\lambda}$, it follows that for an absolute constant $c>0$,
\begin{equation}
 (1-\rho)\left(rT+\sqrt{T/\lambda}\right)\ge c.
 \label{eq:fixed-radius-cohort-mean}
\end{equation}

Now we will show that \eqref{eq:fixed-radius-cohort-mean} also holds when $d<1$. We will construct a lower bound for all the terms in the expression to obtain this lower bound. Using $T\ge T_{\rm queue}$, \eqref{eq:queue-benchmark},
and $d=\lambda T$, we obtain
\begin{equation*}
d=\lambda T
\ge \lambda T_{\rm queue}
=\lambda\bar s+\frac{\lambda^2\E{S^2}}{2(1-\rho)}.
\end{equation*}
Since $\E{S^2}-\bar s^2=\Var(S)\ge0$ and $\rho=\lambda\bar s$, we have $\lambda^2\E{S^2}\ge\lambda^2\bar s^2=\rho^2$.
Consequently,
\begin{equation*}
d\ge\rho+\frac{\rho^2}{2(1-\rho)}.
\end{equation*}
If $1-\rho\le1/3$, then $\rho^2\ge4/9$, so
\begin{equation*}
d\ge\rho+\frac{\rho^2}{2(1-\rho)}
\ge\frac23+\frac{4/9}{2/3}
=\frac43>1.
\end{equation*}
This contradicts $d<1$. Therefore $1-\rho>\frac13$.
 Now we will construct the lower bound for system time terms. Any demand's horizontal coordinate is uniform and independent of the vehicle's position just before arrival, and this horizontal distance gives a lower bound on the distance between two locations. Let $U$ be the arriving demand's horizontal coordinate and $X$ the vehicle's horizontal coordinate just before arrival. Since $U$ is uniform on $[0,1]$ and independent of $X$,
\begin{equation}
\E{|U-X|\mid X=x}
=\int_0^1|u-x|\,\diff u
=\left(x-\frac12\right)^2+\frac14
\ge\frac14.
\end{equation}
Taking expectations over $X$ gives $\E{|U-X|}\ge1/4$. At unit speed, reaching the service disk therefore takes at least $1/4-r$ time units on average. Thus  by $r\leq r_0$, we have $T\ge3/16$, and since $d<1$, we obtain $\sqrt{T/\lambda}=T/\sqrt d>T$. Combining all of these terms gives \eqref{eq:fixed-radius-cohort-mean}, with a constant greater than $c$ in \eqref{eq:fixed-radius-cohort-mean}. Hence \eqref{eq:fixed-radius-cohort-mean} holds also in this case.

By \eqref{eq:fixed-radius-cohort-mean}, at least one of $r(1-\rho)T$ and $(1-\rho)\sqrt{T/\lambda}$ is at least $c/2$. The first gives $T\ge c/[2r(1-\rho)]$ and the second gives $T\ge c^2\lambda/[4(1-\rho)^2]$. Taking the smaller of these bounds proves \eqref{eq:rate-explicit-lower} with $c_0=\min\{c/2,c^2/4\}$.\qed

\begin{proof}[Proof of \Cref{lem:cohort-area}]
Partition the square into $k^2$ equal square cells, with $k=\lfloor8\sqrt d\rfloor$. Since $d\ge1$, we have $7\sqrt d\le k\le8\sqrt d$. During the deterministic interval $[-64T,0]$, each cell receives a Poisson number of arrivals with mean $64d/k^2\ge1$. Thus each cell receives at least one arrival with probability at least $1-e^{-1}$.

Let $J$ count the cells containing the location of an arrival from this interval that has completed service by time zero. Every cell that received arrivals but is omitted from $J$ contains at least one demand still outstanding at time zero. There are at most $N(0)$ such cells, so the expected fraction counted by $J$ satisfies
\[
 \frac{\E{J}}{k^2}\ge 1-e^{-1}-\frac{\E{N(0)}}{k^2}
 \ge 1-e^{-1}-\frac{1}{49}>\frac12.
\]
The vehicle came within $r$ of a demand in each of these $J$ cells. Each cell has diameter $\sqrt{2}/k\le1/(4\sqrt d)$, so it lies entirely within distance $a$ of the realized path, even if the vehicle chose its route after observing the arrivals. The cells have disjoint interiors and area $1/k^2$ each, which gives \eqref{eq:cohort-area}.
\end{proof}

\subsection[Proof of Theorem~\ref*{thm:fixedradius}]{Proof of Theorem~\ref*{thm:fixedradius}}\label{sec:fixedroute}

Stability of the fixed direction and fair coin rules follows from \Cref{prop:stability}. Proofs for forward or reverse direction rules are almost identical. Suppose first that the vehicle uses the forward direction rule. Conditional on the stationary batch size $Q=q$, demands are independent and uniformly distributed. The mean travel distance to each of these arrivals is given by $\E{t(X)}$. Using uniformity and independence, summing over the batch of $q$ demands gives $\E{\sum_{j=1}^Q h_j\mid Q=q}=q\E{t(X)}$. Averaging over $Q$ and substituting into \eqref{eq:rewardbound} gives
\begin{align}
 T_{\rm FR}
 &=\bar s+\frac{1+\rho}{2\lambda}(R_Q-1)+\E{t(X)} \notag\\
 &=T_{\rm queue}+\frac{L_r}{1-\rho}+\E{t(X)}-\frac{L_r}{2},
 \label{eq:fixed-travel-exact}
\end{align}
where the second equality uses \eqref{eq:gatesecond} and \eqref{eq:queue-benchmark}. If the fixed direction is reverse, the vehicle reaches the service point for $x$ after traveling $L_r-t(x)$ rather than $t(x)$, so the only difference is that the mean travel term in \eqref{eq:fixed-travel-exact} becomes $L_r-\E{t(X)}$. For the fair coin rule, substituting mean travel $L_r/2$ into \eqref{eq:rewardbound} gives \eqref{eq:fixedT}.\qed

\subsection[Proof of Theorem~\ref*{thm:scaling}]{Proof of Theorem~\ref*{thm:scaling}}\label{sec:upper}

Fix $0<c<1$ and consider any $c\le\rho<1$ and radius $r>0$. Note that by our model assumptions, $0<m\le M<\infty$ bound the arrival rates across these loads, so $m\le\lambda\le M$. We first give a bound on the travel time on the grid route. We count the segments between successive cell centers, including the return to the depot. There are $k^2$ segments of length $1/k$, so the total length is $k$. Since $k$ is the smallest positive even integer at least $1/(\sqrt{2}r)$, $L_r$ satisfies
\begin{equation}
 \frac{1}{\sqrt{2}\,r}\le L_r=k<\frac{1}{\sqrt{2}\,r}+2.
 \label{eq:uniformcover}
\end{equation}
Therefore, for the grid route policy with a fair coin, \eqref{eq:fixedT} and the length bound \eqref{eq:uniformcover} give
\begin{equation}
 T_{\rm FR}=T_{\rm queue}+\frac{L_r}{1-\rho}
 \le \frac{C_{\rm cov}(1+1/r)}{1-\rho},
\label{eq:frbound}
\end{equation}
 where $C_{\rm cov}$ depends only on $m$ and $K$. The second component of the upper bound is given by an exact location policy. In particular,  we use the fixed batch Traveling Salesperson Problem policy of Bertsimas and van Ryzin~\cite[Section~4.3]{BvR1991}.

Under this policy, each batch consists of $n$ consecutive arrivals, where $n$ is a parameter to be chosen. The vehicle waits at the depot for a batch to be formed, and serves batches in the order they form, using the shortest tour from the depot through the $n$ locations and back to the depot. Since the service is completed at the exact points, this procedure gives a feasible policy for any radius $r\geq0$.

The system time has three components. A demand first waits for its batch to form and then for its batch to depart. After the departure, the remaining component of the system time is the time from departure to service completion. Let $H_n$ denote the length of the tour this demand is assigned to. The last component has a mean of at most $\E{H_n}+n\bar{s}$. Let $W_n$ denote the stationary mean time from the formation of a batch until the vehicle departs to serve it. We will now construct an upper bound on this quantity. The shortest path bound in \cite{Few1955} implies that there exists a positive constant $C$ such that $H_n\le C\sqrt n$. As sums of exponential interarrival times, the intervals between batch formations are independent and Erlang distributed with shape $n$ and rate $\lambda$.
Choose $n=\lceil(2C\lambda/(1-\rho))^2\rceil$ as the batch size parameter. The resulting GI/G/1 queue is stable, since its load satisfies
\[
 \frac{\lambda}{n}\bigl(n\bar s+\E{H_n}\bigr)
 \le \rho+\frac{\lambda C}{\sqrt n}\le \frac{1+\rho}{2}<1.
\]
For $W_n$, \cite[Eq.~(35)]{BvR1991} gives
\[
 W_n\le
 \frac{\lambda\bigl(\lambda^{-2}+\Var(H_n)/n+\Var(S)\bigr)}
      {2\bigl(1-\rho-\lambda\E{H_n}/n\bigr)}
 \le \frac{M\bigl(C^2+(1+K)/m^2\bigr)}{1-\rho},
\]
where the second inequality uses the load bound, $\Var(H_n)/n\le C^2$, and $\Var(S)\le K\bar s^2\le K/m^2$.

For the batch formation waiting component, the demand in arrival position $j$ waits for $n-j$ more arrivals to complete its batch.  Averaging over the $n$ positions gives a mean formation wait of $(n-1)/(2\lambda)$.  Combining these three components gives
\begin{equation}
 T_*\le\frac{n-1}{2\lambda}+W_n+n\bar s+\E{H_n}
 \le\frac{C_{\rm pt}}{(1-\rho)^2},
 \label{eq:pointuniform}
\end{equation}
where $n\le 1+(2CM/(1-\rho))^2$ and $C_{\rm pt}$ depends only on $m,M,K$.

 With these two policies, the upper bound is straightforward to obtain. For $r\ge1-\rho$, use the grid route policy and the upper bound derived in \eqref{eq:frbound}, and for $r<1-\rho$, use the fixed batch bound derived in \eqref{eq:pointuniform}. The optimal mean system time therefore satisfies
\begin{equation*}
 T_*\le\frac{C_{m,M,K}}{1-\rho}
 \left(1+\min\left\{\frac1r,\frac1{1-\rho}\right\}\right),
\end{equation*}
where $C_{m,M,K}$ can be taken as $\max\{C_{\rm cov},C_{\rm pt}\}$.  This establishes the desired order of the upper bound.

We now prove the matching lower bound. \Cref{thm:fixed-radius-lower} gives $T\ge T_{\rm queue}$ for every admissible stable policy. Since $\lambda\E{S^2}\ge\lambda\bar s^2=\rho^2/\lambda\ge c^2/M$, \eqref{eq:queue-benchmark} gives
\[
 T\ge\frac{c^2}{2M(1-\rho)}.
\]
There are two cases to consider. For $r\le r_0$, using the lower bound on $\lambda$ in \eqref{eq:rate-explicit-lower} also gives
\[
 T\ge \frac{c_0\min\{1,m\}}{1-\rho}\min\left\{\frac1r,\frac1{1-\rho}\right\}.
\]
For $r>r_0$, we have
\begin{equation}
1+\min\left\{\frac1r,\frac1{1-\rho}\right\}
\le 1+\frac1{r_0}.
\end{equation}
Then it holds that
\begin{equation}
T\ge \frac{c^2}{2M(1-\rho)}
\ge \frac{c^2}{2M(1+1/r_0)(1-\rho)}
\left(1+\min\left\{\frac1r,\frac1{1-\rho}\right\}\right).
\end{equation}
This establishes the required lower bound order for $r>r_0$. Since the policy was selected arbitrarily, the lower bounds above therefore give the desired lower bound for $T_*$.

At $r=0$, the fixed batch policy gives the upper bound \eqref{eq:pointuniform}. Every policy feasible at radius zero is also feasible at $r=1-\rho$, so the lower bound just proved at $r=1-\rho$ gives the matching order $(1-\rho)^{-2}$. \qed

\subsection[Proof of Theorem~\ref*{thm:fixed-route-scaling}]{Proof of Theorem~\ref*{thm:fixed-route-scaling}}\label{sec:proof-fixed-route-scaling}

To determine whether a fixed-route policy is order optimal, we compare its mean system time with $T_*$. Since there is an exact formula for the mean system time of the fair coin rule, we first show that any stable direction rule on the same route has asymptotically the same mean system time (their ratio tends to one as $\rho\rightarrow1$).

Consider such a direction rule and its stationary regime. By \eqref{eq:gatesecond}, $R_Q$ is the same under this rule and the fair coin rule, and only the travel term in \eqref{eq:rewardbound} can differ. Since $0\le h_j\le L_r$ for every $j$, the travel term lies between zero and $L_r$, and it equals $L_r/2$ under the fair coin rule. The two mean system times therefore differ by at most $L_r/2$. The fair coin mean in \eqref{eq:fixedT} is at least $L_r/(1-\rho)$, so the relative difference is at most $(1-\rho)/2$ for every route. The relative difference vanishes as $\rho\rightarrow1$. It therefore suffices to consider the fair coin rule for the rest of the proof.

When $\alpha>0$, $r=(1-\rho)^\alpha$ tends to zero. The $r$-neighborhood of the route covers the unit square, so \eqref{eq:path-area} gives $1\le 2rL_r+\pi r^2$. For $r$ small enough that $\pi r^2\le1/2$, this implies $L_r\ge1/(4r)$. Together with the assumed upper bound $L_r=O(1/r)$, we get $L_r=\Theta(1/r)$.

By \eqref{eq:queue-benchmark} and \eqref{eq:fixedT}, $(1-\rho)T_{\rm FR}=L_r+\Theta(1)$. This is of order $1/r$ for $\alpha>0$ and of order one at $\alpha=0$, since then $r=1$ and $L_r=O(1)$. Thus
\[
 T_{\text{FR}}(\alpha)=\Theta\!\left(\frac{1}{r(1-\rho)}\right)
 =\Theta\!\left((1-\rho)^{-(1+\alpha)}\right).
\]
Comparing with \eqref{eq:alphascaling} gives the stated range of order optimality. If $r/(1-\rho)\to0$, every covering route has $L_r\ge1/(4r)$ eventually. By \eqref{eq:fixedT}, $T_{\rm FR}\ge 1/[4r(1-\rho)]$, while \Cref{thm:scaling} gives $T_*=\Theta((1-\rho)^{-2})$. Thus $T_{\rm FR}/T_*\to\infty$.\qed

\fi

\section{Conclusion}\label{sec:conclusion}

Service location flexibility lets a vehicle serve nearby customers at a shared stop, reducing travel and the wait for other customers. We show that fixed-route policies can attain the optimal order of mean system time in heavy traffic when the service radius is of order $1-\rho$ or larger. The grid route policy is one representative policy in this class. It follows a closed route through the centers of square cells and serves each demand at the center of its cell. These service points can be lockers or other pickup points in parcel delivery, or regional service centers in a mobile service program. The carrier chooses the service points and route together before demands arrive.

Moreover, for fixed-direction and independent fair-coin policies, we characterize precisely how mean system time depends on the chosen route at every stable load, including intermediate loads. In particular, if the vehicle tosses an independent fair coin at each departure from the depot to choose the forward or reverse direction, mean system time increases linearly with route length when the arrival process and service-time distribution are fixed.

We derived these savings assuming that all customers accept the same pickup distance. Whether customers would accept pickup distances that yield substantial savings in a particular delivery system remains an empirical question. In practice, a carrier can offer a discount for collection at a pickup point~\cite{Galiullina2024}. The customer who accepts incurs the extra travel, while a shorter route can also reduce the wait for later orders. With a fixed budget on the sum of the service radii, would the carrier save more time by asking many customers to walk a little or a few to walk farther? A related allocation question arises in spatial matching~\cite{AmeenYu2026}. A carrier might also offer different pickup distances in dense and sparse neighborhoods, where the distances between orders differ.

\let\originalthebibliography\thebibliography
\renewcommand{\thebibliography}[1]{%
    \originalthebibliography{#1}%
    \phantomsection
    \addcontentsline{toc}{section}{References}%
}
\bibliographystyle{alpha}
\bibliography{bibliography}

@article{BartholdiPlatzman1988,
  author = {Bartholdi, III, John J. and Platzman, Loren K.},
  title = {Heuristics Based on Spacefilling Curves for Combinatorial Problems in {Euclidean} Space},
  journal = {Management Science},
  volume = {34},
  number = {3},
  pages = {291--305},
  year = {1988},
  doi = {10.1287/mnsc.34.3.291}
}

@article{BvR1991,
author = {Bertsimas, Dimitris J. and van Ryzin, Garrett},
title = {A Stochastic and Dynamic Vehicle Routing Problem in the {Euclidean} Plane},
journal = {Operations Research},
volume = {39},
number = {4},
pages = {601-615},
year = {1991},
doi = {10.1287/opre.39.4.601},
URL = {https://doi.org/10.1287/opre.39.4.601},
eprint = { https://doi.org/10.1287/opre.39.4.601}}

@article{Pavone2011,
  author={Pavone, Marco and Frazzoli, Emilio and Bullo, Francesco},
  journal={IEEE Transactions on Automatic Control}, 
  title={Adaptive and Distributed Algorithms for Vehicle Routing in a Stochastic and Dynamic Environment}, 
  year={2011},
  volume={56},
  number={6},
  pages={1259-1274},
  doi={10.1109/TAC.2010.2092850}}

@techreport{Carlsson,
  author       = {Carlsson, John Gunnar and Meng, Xiangfei},
  title        = {Quantifying the Impact of Next-Generation Modes of Delivery},
  institution  = {METRANS Transportation Center},
  number       = {15-14},
  year         = {2016},
  url          = {https://rosap.ntl.bts.gov/view/dot/67577}
}

@article{CarlssonEtAl2016,
  author  = {Carlsson, John Gunnar and Behroozi, Mehdi and Devulapalli, Raghuveer and Meng, Xiangfei},
  title   = {Household-Level Economies of Scale in Transportation},
  journal = {Operations Research},
  volume  = {64},
  number  = {6},
  pages   = {1372--1387},
  year    = {2016},
  doi     = {10.1287/opre.2016.1533},
  url     = {https://doi.org/10.1287/opre.2016.1533}
}

@inproceedings{Bandyapadhyay2025,
  author    = {Bandyapadhyay, Sayan and Clinch, Katie and Lochet, William and Lokshtanov, Daniel and Saurabh, Saket and Xue, Jie},
  title     = {{PTASes} for {Euclidean TSP} with Unit Disk and Unit Square Neighborhoods},
  booktitle = {Proceedings of the 2025 Annual ACM-SIAM Symposium on Discrete Algorithms},
  pages     = {2326--2356},
  publisher = {SIAM},
  year      = {2025},
  doi       = {10.1137/1.9781611978322.78}
}

@incollection{Steele1990,
  author    = {Steele, J. Michael},
  title     = {Seedlings in the Theory of Shortest Paths},
  booktitle = {Disorder in Physical Systems:
               A Volume in Honour of John M. Hammersley},
  editor    = {Grimmett, Geoffrey R. and Welsh, Dominic J. A.},
  pages     = {277--306},
  publisher = {Oxford University Press},
  address   = {Oxford},
  year      = {1990},
  url       = {https://www.statslab.cam.ac.uk/~grg1000/books/hammfest/17-jms.pdf}
}

@article{Beardwood_Halton_Hammersley_1959, title={The shortest path through many points}, volume={55}, DOI={10.1017/S0305004100034095}, number={4}, journal={Proceedings of the Cambridge Philosophical Society}, author={Beardwood, Jillian and Halton, J. H. and Hammersley, J. M.}, year={1959}, pages={299–327}}

@phdthesis{xu1994optimal,
  title={Optimal policies for stochastic and dynamic vehicle routing problems},
  author={Xu, Haiping},
  year={1995},
  school={Massachusetts Institute of Technology},
  url={https://hdl.handle.net/1721.1/11748}
}

@article{Celik2012,
  title = {Controlled mobility in stochastic and dynamic wireless networks},
  volume = {72},
  ISSN = {1572-9443},
  url = {http://dx.doi.org/10.1007/s11134-012-9313-y},
  DOI = {10.1007/s11134-012-9313-y},
  number = {3-4},
  journal = {Queueing Systems},
  publisher = {Springer Science and Business Media LLC},
  author = {{\c C}elik, G{\"u}ner D. and Modiano, Eytan H.},
  year = {2012},
  month = dec,
  pages = {251–277}
}

@inproceedings{Fekete,
author = {Sándor P. Fekete and Dominik Krupke and Michael Perk and Christian Rieck and Christian Scheffer},
title = {A Closer Cut: Computing Near-Optimal Lawn Mowing Tours},
booktitle = {2023 Proceedings of the Symposium on Algorithm Engineering and Experiments (ALENEX)},
publisher = {SIAM},
year={2023},
pages = {1-14},
doi = {10.1137/1.9781611977561.ch1},
URL = {https://epubs.siam.org/doi/abs/10.1137/1.9781611977561.ch1},
eprint = {https://epubs.siam.org/doi/pdf/10.1137/1.9781611977561.ch1}}

@article{ARKIN200025,
title = {Approximation algorithms for lawn mowing and milling},
journal = {Computational Geometry},
volume = {17},
number = {1--2},
pages = {25--50},
year = {2000},
issn = {0925-7721},
doi = {10.1016/S0925-7721(00)00015-8},
url = {https://www.sciencedirect.com/science/article/pii/S0925772100000158},
author = {Esther M. Arkin and Sándor P. Fekete and Joseph S. B. Mitchell},
}

@article{Bertsimas1993multi,
 ISSN = {0030364X, 15265463},
 URL = {http://www.jstor.org/stable/171944},
 author = {Dimitris J. Bertsimas and Garrett van Ryzin},
 journal = {Operations Research},
 number = {1},
 pages = {60--76},
 publisher = {INFORMS},
 title = {Stochastic and Dynamic Vehicle Routing in the {Euclidean} Plane with Multiple Capacitated Vehicles},
 doi = {10.1287/opre.41.1.60},
 urldate = {2026-08-05},
 volume = {41},
 year = {1993}
}

@article{Bertsimas1993,
  title = {Stochastic and dynamic vehicle routing with general demand and interarrival time distributions},
  volume = {25},
  ISSN = {1475-6064},
  url = {http://dx.doi.org/10.2307/1427801},
  DOI = {10.2307/1427801},
  number = {4},
  journal = {Advances in Applied Probability},
  publisher = {Cambridge University Press (CUP)},
  author = {Bertsimas,  Dimitris J. and van Ryzin,  Garrett},
  year = {1993},
  month = Dec,
  pages = {947–978}
}

@misc{dhlecommerce2026ooh,
  author       = {{DHL eCommerce}},
  title        = {Out-of-Home Delivery \& Returns -- {Europe} Leads the Way!},
  year         = {2026},
  url          = {https://www.dhl.com/global-en/microsites/ec/ecommerce-insights/insights/e-commerce-logistics/out-of-home-trends.html},
  howpublished = {\href{https://www.dhl.com/global-en/microsites/ec/ecommerce-insights/insights/e-commerce-logistics/out-of-home-trends.html}{DHL eCommerce Insights}},
  note         = {Accessed November 9, 2026}
}

@article{Wang2025,
author = {Wang, Quanmeng and Lyu, Guodong and He, Long and Teo, Chung Piaw},
title = {Does the Locker Alliance Network Improve Last Mile Delivery Efficiency?},
journal = {Management Science},
volume = {72},
number = {5},
pages = {4124-4143},
year = {2026},
doi = {10.1287/mnsc.2024.04567},

URL = { 
    
        https://doi.org/10.1287/mnsc.2024.04567
    
    

},
eprint = {    https://doi.org/10.1287/mnsc.2024.04567}
}

@article{Galiullina2024,
author = {Galiullina, Albina and Mutlu, Nevin and Kinable, Joris and Van Woensel, Tom},
title = {Demand Steering in a Last-Mile Delivery Problem with Home and Pickup Point Delivery Options},
journal = {Transportation Science},
volume = {58},
number = {2},
pages = {454-473},
year = {2024},
doi = {10.1287/trsc.2023.0287},

URL = { 
    
        https://doi.org/10.1287/trsc.2023.0287
    
    

},
eprint = { 
    
        https://doi.org/10.1287/trsc.2023.0287

}}

@article{Glaeser2019,
author = {Glaeser, Chloe Kim and Fisher, Marshall and Su, Xuanming},
title = {Optimal Retail Location: Empirical Methodology and Application to Practice},
journal = {Manufacturing \& Service Operations Management},
volume = {21},
number = {1},
pages = {86-102},
year = {2019},
doi = {10.1287/msom.2018.0759},

URL = { 
    
        https://doi.org/10.1287/msom.2018.0759
    
    

},
eprint = { 
    
        https://doi.org/10.1287/msom.2018.0759
    
    

}}

@article{Janinhoff2024,
title = {Out-of-home delivery in last-mile logistics: A review},
journal = {Computers \& Operations Research},
volume = {168},
pages = {106686},
year = {2024},
issn = {0305-0548},
doi = {10.1016/j.cor.2024.106686},
url = {https://www.sciencedirect.com/science/article/pii/S0305054824001588},
author = {Lukas Janinhoff and Robert Klein and Daniela Sailer and Jim Morten Schoppa},
}

@article{Dumitrescu2003,
title = {Approximation algorithms for {TSP} with neighborhoods in the plane},
journal = {Journal of Algorithms},
volume = {48},
number = {1},
pages = {135-159},
year = {2003},
issn = {0196-6774},
doi = {10.1016/S0196-6774(03)00047-6},
url = {https://www.sciencedirect.com/science/article/pii/S0196677403000476},
author = {Adrian Dumitrescu and Joseph S. B. Mitchell},
}

@article{Tilk2021,
  title = {The last-mile vehicle routing problem with delivery options},
  volume = {43},
  ISSN = {1436-6304},
  url = {http://dx.doi.org/10.1007/s00291-021-00633-0},
  DOI = {10.1007/s00291-021-00633-0},
  number = {4},
  journal = {OR Spectrum},
  publisher = {Springer Science and Business Media LLC},
  author = {Tilk,  Christian and Olkis,  Katharina and Irnich,  Stefan},
  year = {2021},
  month = dec,
  pages = {877–904}
}

@article{Dumez2021,
title = {A large neighborhood search approach to the vehicle routing problem with delivery options},
journal = {Transportation Research Part B: Methodological},
volume = {144},
pages = {103-132},
year = {2021},
issn = {0191-2615},
doi = {10.1016/j.trb.2020.11.012},
url = {https://www.sciencedirect.com/science/article/pii/S0191261520304422},
author = {Dorian Dumez and Fabien Lehuédé and Olivier Péton},
}

@article{Behdani2014,
author = {Behdani, Behnam and Smith, J. Cole},
title = {An Integer-Programming-Based Approach to the Close-Enough Traveling Salesman Problem},
journal = {INFORMS Journal on Computing},
volume = {26},
number = {3},
pages = {415-432},
year = {2014},
doi = {10.1287/ijoc.2013.0574},

URL = { 
    
        https://doi.org/10.1287/ijoc.2013.0574
    
    

},
eprint = { 
    
        https://doi.org/10.1287/ijoc.2013.0574
}}

@article{Psaraftis1995,
  title = {Dynamic vehicle routing: Status and prospects},
  volume = {61},
  ISSN = {1572-9338},
  url = {http://dx.doi.org/10.1007/BF02098286},
  DOI = {10.1007/bf02098286},
  number = {1},
  journal = {Annals of Operations Research},
  publisher = {Springer Science and Business Media LLC},
  author = {Psaraftis,  Harilaos N.},
  year = {1995},
  month = Dec,
  pages = {143–164}
}

@article{Psaraftis1980,
author = {Psaraftis, Harilaos N.},
title = {A Dynamic Programming Solution to the Single Vehicle Many-to-Many Immediate Request Dial-a-Ride Problem},
journal = {Transportation Science},
volume = {14},
number = {2},
pages = {130-154},
year = {1980},
doi = {10.1287/trsc.14.2.130},

URL = { 
    
        https://doi.org/10.1287/trsc.14.2.130
    
    

},
eprint = { 
    
        https://doi.org/10.1287/trsc.14.2.130
    
    

}}

@article{Berman1985,
  title = {Optimal Server Location on a Network Operating as an {M/G/1} Queue},
  volume = {33},
  ISSN = {1526-5463},
  url = {http://dx.doi.org/10.1287/opre.33.4.746},
  DOI = {10.1287/opre.33.4.746},
  number = {4},
  journal = {Operations Research},
  publisher = {Institute for Operations Research and the Management Sciences (INFORMS)},
  author = {Berman,  Oded and Larson,  Richard C. and Chiu,  Samuel S.},
  year = {1985},
  month = Aug,
  pages = {746–771}
}

@article{KavithaAltman2012,
  author = {Kavitha, Veeraruna and Altman, Eitan},
  title = {Continuous polling models and application to ferry assisted {WLAN}},
  journal = {Annals of Operations Research},
  year = {2012},
  volume = {198},
  number = {1},
  pages = {185--218},
  doi = {10.1007/s10479-011-0960-x}
}

@article{Altman1994,
  title = {Queueing in space},
  volume = {26},
  ISSN = {1475-6064},
  url = {http://dx.doi.org/10.2307/1427906},
  DOI = {10.2307/1427906},
  number = {4},
  journal = {Advances in Applied Probability},
  publisher = {Cambridge University Press (CUP)},
  author = {Altman,  Eitan and Levy,  Hanoch},
  year = {1994},
  month = Dec,
  pages = {1095–1116}
}

@misc{Freund2024,
  doi = {10.48550/ARXIV.2404.04709},
  url = {https://arxiv.org/abs/2404.04709},
  author = {Freund,  Daniel and Martin,  Sébastien and Zhao,  Jiayu Kamessi},
  title = {Two-Sided Flexibility in Platforms},
  publisher = {arXiv},
  howpublished = {arXiv preprint \href{https://arxiv.org/abs/2404.04709}{arXiv:2404.04709}},
  year = {2024},
  copyright = {Creative Commons Attribution 4.0 International}
}

@misc{AmeenYu2026,
  author = {Ameen, Taha and Yu, Sophie H.},
  title = {Distance flexibility in spatial matching: the value of concentration},
  year = {2026},
  doi = {10.48550/arXiv.2609.36361},
  url = {https://arxiv.org/abs/2609.36361},
  howpublished = {arXiv preprint \href{https://arxiv.org/abs/2609.36361}{arXiv:2609.36361}}
}

@misc{Uniformity2026,
  doi = {10.48550/ARXIV.2601.13426},
  url = {https://arxiv.org/abs/2601.13426},
  author = {Ameen,  Taha and Sentenac,  Flore and Yu,  Sophie H.},
  title = {A uniformity principle for spatial matching},
  publisher = {arXiv},
  howpublished = {arXiv preprint \href{https://arxiv.org/abs/2601.13426}{arXiv:2601.13426}},
  year = {2026},
  copyright = {arXiv.org perpetual,  non-exclusive license}
}

@misc{Flexibility2026,
  doi = {10.48550/ARXIV.2604.02295},
  url = {https://arxiv.org/abs/2604.02295},
  author = {Ameen,  Taha and Sentenac,  Flore and Yu,  Sophie H.},
  title = {Flexibility allocation in random bipartite matching markets: exact matching rates and dominance regimes},
  publisher = {arXiv},
  howpublished = {arXiv preprint \href{https://arxiv.org/abs/2604.02295}{arXiv:2604.02295}},
  year = {2026},
  copyright = {arXiv.org perpetual,  non-exclusive license}
}

@article{Krushinsky2021,
  title = {Location flexibility in parcel delivery operations: framework and empirical analysis},
  volume = {33},
  ISSN = {1471-6798},
  url = {http://dx.doi.org/10.1093/imaman/dpab004},
  DOI = {10.1093/imaman/dpab004},
  number = {4},
  journal = {IMA Journal of Management Mathematics},
  publisher = {Oxford University Press (OUP)},
  author = {Krushinsky,  Dmitry and Guo,  Xuezhen and Claassen,  G. D. H.},
  year = {2022},
  month = oct,
  pages = {563–581}
}

@article{Pillac2013,
  title = {A review of dynamic vehicle routing problems},
  volume = {225},
  ISSN = {0377-2217},
  url = {http://dx.doi.org/10.1016/j.ejor.2012.08.015},
  DOI = {10.1016/j.ejor.2012.08.015},
  number = {1},
  journal = {European Journal of Operational Research},
  publisher = {Elsevier BV},
  author = {Pillac,  Victor and Gendreau,  Michel and Guéret,  Christelle and Medaglia,  Andrés L.},
  year = {2013},
  month = Feb,
  pages = {1–11}
}

@book{MeynTweedie1993,
  author = {Meyn, Sean P. and Tweedie, Richard L.},
  title = {Markov Chains and Stochastic Stability},
  publisher = {Springer-Verlag},
  address = {London},
  year = {1993},
  doi = {10.1007/978-1-4471-3267-7},
  url = {https://probability.ca/MT/}
}

@book{Thorisson2000,
  author = {Thorisson, Hermann},
  title = {Coupling, Stationarity, and Regeneration},
  publisher = {Springer-Verlag},
  address = {New York},
  year = {2000},
  doi = {10.1007/978-1-4612-1236-2}
}

@book{BertsekasGallager1992,
  author = {Bertsekas, Dimitri P. and Gallager, Robert G.},
  title = {Data Networks},
  edition = {Second},
  publisher = {Prentice Hall},
  year = {1992},
  url = {https://web.mit.edu/dimitrib/www/datanets.html}
}

@article{Kroese1997,
  author = {Kroese, Dirk P.},
  title = {Heavy traffic analysis for continuous polling models},
  journal = {Journal of Applied Probability},
  year = {1997},
  volume = {34},
  number = {3},
  pages = {720--732},
  doi = {10.2307/3215097}
}

@article{BoxmaLevyYechiali1992,
  author = {Boxma, Onno J. and Levy, Hanoch and Yechiali, Uri},
  title = {Cyclic reservation schemes for efficient operation of multiple-queue single-server systems},
  journal = {Annals of Operations Research},
  year = {1992},
  volume = {35},
  pages = {187--208},
  doi = {10.1007/BF02188704}
}

@article{Few1955,
  author = {Few, L.},
  title = {The shortest path and the shortest road through $n$ points},
  journal = {Mathematika},
  year = {1955},
  volume = {2},
  number = {2},
  pages = {141--144},
  doi = {10.1112/S0025579300000784}
}

@article{UlmerStreng2019,
  author = {Ulmer, Marlin W. and Streng, Sebastian},
  title = {Same-day delivery with pickup stations and autonomous vehicles},
  journal = {Computers \& Operations Research},
  year = {2019},
  volume = {108},
  pages = {1--19},
  doi = {10.1016/j.cor.2019.03.017}
}

@article{FuhrmannCooper1985Continuum,
  author = {Fuhrmann, S. W. and Cooper, R. B.},
  title = {Application of decomposition principle in {M/G/1} vacation model to two continuum cyclic queueing models---especially token-ring {LAN}s},
  journal = {AT\&T Technical Journal},
  year = {1985},
  volume = {64},
  number = {5},
  pages = {1091--1099},
  doi = {10.1002/j.1538-7305.1985.tb00457.x}
}

\end{document}